\newcommand{\xBc}{\langle}
\newcommand{\xBe}{\rangle}

\newcommand{\xbP}{\Pi}

\newcommand{\xbS}{\Sigma}

\newcommand{\xba}{\alpha}
\newcommand{\xbb}{\beta}

\newcommand{\xbe}{\in}
\newcommand{\xbf}{\phi}

\newcommand{\xbm}{\mu}

\newcommand{\xbo}{\omega}

\newcommand{\xbq}{\psi}

\newcommand{\xbs}{\sigma}
\newcommand{\xbt}{\tau}

\newcommand{\xCB}{A}

\newcommand{\xCK}{\times}

\newcommand{\xCN}{\neg}
\newcommand{\xCO}{ }
\newcommand{\xCQ}{\emptyset}

\newcommand{\xCf}{\hspace{0.1em}}

\newcommand{\xcA}{\forall}

\newcommand{\xcC}{\not\subseteq}

\newcommand{\xcE}{\exists}

\newcommand{\xcM}{\not\models}

\newcommand{\xcS}{\bigcap}

\newcommand{\xcV}{\bigcup}

\newcommand{\xcc}{\subseteq}

\newcommand{\xce}{\not\in}

\newcommand{\xch}{\Rightarrow}

\newcommand{\xck}{\leq}
\newcommand{\xcl}{\vdash}
\newcommand{\xcm}{\models}
\newcommand{\xcn}{\hspace{0.2em}\sim\hspace{-0.9em}\mid\hspace{0.58em}}

\newcommand{\xcp}{\rightarrow}

\newcommand{\xcr}{\leftrightarrow}
\newcommand{\xcs}{\cap}
\newcommand{\xcu}{\wedge}
\newcommand{\xcv}{\cup}

\newcommand{\xcz}{\Box}

\newcommand{\xDH}{\item }

\newcommand{\xda}{{\cal A}}
\newcommand{\xdb}{{\cal B}}

\newcommand{\xdl}{{\cal L}}

\newcommand{\xdx}{{\cal X}}
\newcommand{\xdy}{{\cal Y}}

\newcommand{\xEI}{\begin{itemize}}
\newcommand{\xEJ}{\end{itemize}}

\newcommand{\xEd}{\neq}
\newcommand{\xEh}{\begin{enumerate}}
\newcommand{\xEj}{\end{enumerate}}

\newcommand{\xeb}{\prec}

\newcommand{\xex}{\lceil}

\newcommand{\xfA}{\mid}

\newcommand{\Xl}{\ldots}

\newcommand{\bl}{\begin{lemma} \rm}
\newcommand{\el}{\end{lemma}}
\newcommand{\br}{\begin{remark} \rm}
\newcommand{\er}{\end{remark}}
\newcommand{\be}{\begin{example} \rm}
\newcommand{\ee}{\end{example}}
\newcommand{\bco}{\begin{corollary} \rm}
\newcommand{\eco}{\end{corollary}}
\newcommand{\bc}{\begin{claim} \rm}
\newcommand{\ec}{\end{claim}}
\newcommand{\bfa}{\begin{fact} \rm}
\newcommand{\efa}{\end{fact}}
\newcommand{\bp}{\begin{proposition} \rm}
\newcommand{\ep}{\end{proposition}}
\newcommand{\bd}{\begin{definition} \rm}
\newcommand{\ed}{\end{definition}}
\newcommand{\bcs}{\begin{construction} \rm}
\newcommand{\ecs}{\end{construction}}
\newcommand{\bcd}{\begin{condition} \rm}
\newcommand{\ecd}{\end{condition}}
\newcommand{\bt}{\begin{theorem} \rm}
\newcommand{\et}{\end{theorem}}
\newcommand{\bn}{\begin{notation} \rm}
\newcommand{\en}{\end{notation}}
\newcommand{\bfi}{\begin{bild} \rm}
\newcommand{\efi}{\end{bild}}
\newcommand{\bsta}{\begin{statement} \rm}
\newcommand{\esta}{\end{statement}}
\newcommand{\bcom}{\begin{comment} \rm}
\newcommand{\ecom}{\end{comment}}
\newcommand{\bdia}{\begin{diagram} \rm}
\newcommand{\edia}{\end{diagram}}

\newcommand{\bfc}{\begin{figure}[htb] \begin{center}}
\newcommand{\efc}{\end{center} \end{figure}}

\documentclass{article}

\usepackage{amssymb,latexsym,epic,eepic,rotating}

\title{
Independence - revision and defaults
\thanks{
Paper 329
}
}

\author{Dov M Gabbay
\thanks{
Dov.Gabbay@kcl.ac.uk, www.dcs.kcl.ac.uk/staff/dg
} \\
King's College, London
\thanks{
Department of Computer Science, King's College London, Strand,
London WC2R 2LS, UK
} \\
and \\
Bar-Ilan University, Israel \\ \\
Karl Schlechta
\thanks{
ks@cmi.univ-mrs.fr, karl.schlechta@web.de, http://www.cmi.univ-mrs.fr/ $\sim$ ks
} \\
Laboratoire d'Informatique Fondamentale de Marseille
\thanks{
UMR 6166, CNRS and Universit\'{e} de Provence,
Address: CMI, 39, rue Joliot-Curie, F-13453 Marseille Cedex 13, France
}
}

\begin{document}

\newtheorem{lemma}{Lemma}[section]
\newtheorem{theorem}[lemma]{Theorem}
\newtheorem{proposition}[lemma]{Proposition}
\newtheorem{corollary}[lemma]{Corollary}
\newtheorem{claim}[lemma]{Claim}
\newtheorem{fact}[lemma]{Fact}
\newtheorem{remark}[lemma]{Remark}
\newtheorem{definition}{Definition}[section]
\newtheorem{construction}{Construction}[section]
\newtheorem{condition}{Condition}[section]
\newtheorem{example}{Example}[section]
\newtheorem{notation}{Notation}[section]
\newtheorem{bild}{Figure}[section]
\newtheorem{comment}{Comment}[section]
\newtheorem{statement}{Statement}[section]
\newtheorem{diagram}{Diagram}[section]

\renewcommand{\labelenumi}
  {(\arabic{enumi})}
\renewcommand{\labelenumii}
  {(\arabic{enumi}.\arabic{enumii})}
\renewcommand{\labelenumiii}
  {(\arabic{enumi}.\arabic{enumii}.\arabic{enumiii})}
\renewcommand{\labelenumiv}
  {(\arabic{enumi}.\arabic{enumii}.\arabic{enumiii}.\arabic{enumiv})}

\maketitle

\setcounter{secnumdepth}{3}
\setcounter{tocdepth}{3}

\begin{abstract}

We investigate different aspects of independence here, in the context
of theory revision, generalizing slightly work by Chopra, Parikh, and
Rodrigues, and in the context of preferential reasoning.

\end{abstract}

\tableofcontents

%
%
%


\section{
Introduction
}

We give some results on

 \xEh

 \xDH

Theory revision:

Parikh and co-authors (see  \cite{CP00}), and, independently,
Rodrigues
(see  \cite{Rod97}), have
investigated a notion of logical independence, based
on the sharing of essential propositional variables. We do a semantical
analogue here. What Parikh et al. call splitting on the logical level, we
call
factorization (on the semantical level).

A comparison of the work of Parikh and Rodrigues can be
found in  \cite{Mak09}.

Note that many of our results are
valid for arbitrary products, not only for classical model sets.

We go very slightly beyond Parikh's work, as a matter of fact our
generalization
is alreday contained in the Axiom P2g, due to K.Georgatos (see
 \cite{CP00}):

$ \xCf (P2g)$ If $T$ is split between $ \xdl_{1}$ and $ \xdl_{2},$ $ \xba
,$ $ \xbb $ are in $ \xdl_{1}$ and $ \xdl_{2}$
respectively, then $T* \xba *b=T* \xbb * \xba =T*(\xba \xcu \xbb).$

On the other hand, we stay below Parikh's work and do $ \xCf not$
investigate
partial overlap (see $ \xdb -$structures model
in  \cite{CP00}).

We claim no originality of the basic ideas, just our proofs and perhaps an
example might be new - but they are always elementary and very easy.

 \xDH

Preferential reasoning:

We shortly discuss preferential structures which have properties
of defaults in the fact that they permit to treat sub-ideal information.
Usually, we have only the ideal case, where
all ``normal'' information holds, and the classical case.
(Reiter) defaults, but also e.g., inheritance systems, permit to
satisfy only some, but not necessarily all, default rules, and
are thus more flexible. We show how to construct preferential
structures with the same properties.

 \xEj
\subsection{
The situation in the case of theory revision
}

We work here with arbitrary, non-empty products. Intuitively, $ \xdy $ is
the set
of models for the propositional variable set $U.$ We assume the Axiom of
Choice.

\bd

$\hspace{0.01em}$


\label{Definition Factor-1.1}

Let $U$ be an index set, $ \xdy = \xbP \{Y_{k}:k \xbe U\},$ let all $Y_{k}
\xEd \xCQ,$ and $ \xdx \xcc \xdy.$ Thus, $ \xbs \xbe \xdx $
is a function from $U$ to $ \xcV \{Y_{k}:k \xbe U\}$ s.t. $ \xbs (k) \xbe
Y_{k}.$
We then note $X_{k}:=\{y \xbe Y_{k}: \xcE \xbs \xbe \xdx. \xbs (k)=y\}.$

If $U' \xcc U,$ then $ \xbs \xex U' $ will be the restriction of $ \xbs $
to $U',$ and
$ \xdx \xex U':=\{ \xbs \xex U': \xbs \xbe \xdx \}.$

If $ \xda:=\{A_{i}:i \xbe I\}$ is a partition of $U,$ $U' \xcc U,$ then $
\xda \xex U':=\{A_{i} \xcs U' \xEd \xCQ:i \xbe I\}.$

Let $ \xda:=\{A_{i}:i \xbe I\},$ $ \xdb:=\{B_{j}:j \xbe J\}$ both be
partitions of $U,$ then $ \xda $ is called a
refinement of $ \xdb $ iff for all $i \xbe I$ there is $j \xbe J$ s.t.
$A_{i} \xcc B_{j}.$

A partition $ \xda $ of $U$ will be called a factorization of $ \xdx $ iff
$ \xdx =\{ \xbs \xbe \xdy: \xcA i \xbe I(\xbs \xex A_{i} \xbe \xdx \xex
A_{i})\},$ we will also sometimes say for clarity that
$ \xda $ is a partition of $ \xdx $ over $U.$

We will adhere to above notations throughout these pages.

If $ \xdx $ is as above, $U' \xcc U,$ and $ \xbs \xbe \xdx \xex U',$ then
there is obviously some (usually
not unique) $ \xbt \xbe \xdx $ s.t. $ \xbt \xex U' = \xbs.$ This trivial
fact will be used repeatedly in
the following pages. We will denote by $ \xbs^{+}$ some such $ \xbt $ -
context will
tell which are the $U' $ and $U.$ (To be more definite, we may take the
first such $ \xbt $
in some arbitrary enumeration of $ \xdx.)$

Given a propositional language $ \xdl,$ $v(\xdl)$ will be the set of
its
propositional variables, and $v(\xbf)$ the set of variables occuring in
$ \xbf.$
A model set $C$ is called definable iff there is a theory $T$ s.t.
$C=M(T)$ - the
set of models of $T.$
\subsection{Organization}

\ed

We treat here the following:

 \xEh

 \xDH We give a purely algebraic description of factorization.

This is the algebraic analogue of work by Parikh and co-authors, and we
claim almost no originality, perhaps with the exception of an example and
the remark on language independence.

 \xDH We generalize slightly the Parikh approach so it can be described
as commuting with decomposition into sublanguages and addition.

We show by a trivial argument that this corresponds to a generalized
Hamming
distance between models.

 \xDH We go beyond Rational Monotony and show how to construct a
preferential
structure from a set of (normal) defaults. Thus, we give an independent
semantics to normal defaults, translating their usual treatment into
a homogenous construction of the preferential structure.

In the general case, this gives nothing new, as any preferential structure
can be constructed this way. (We consider the one-copy case only.)
Most of the time, it will result in a special structure which
automatically
takes into account the specificity criterion to resolve conflicts. The
essential idea is to take a modified Hamming distance on the set of
satisfied defaults, modified as we do not count the defaults, but look at
them as sets, together with the subset relation.

We also show that our approach can be seen as a revision of the ideal,
perhaps non-existant case, or as an approach to this ideal case as the
limit. Of course, when the ideal case is consistent, then this will be
our result.

 \xDH Independence in the case of TR is treated by looking at
``independent'' parts ``independently'',
and later summing up. In the case of defaults, we treat the defaults
independently, just as in the Reiter approach, but also
``inside'' the model sets, we treat subsets just as we the sets
themselves, resulting in a partial kind of rankedness (by default).

 \xDH

We conclude by giving a simple informal argument why the TR situation is
more
complicated than the default situation.

 \xEj
\section{
Factorisation
}

\bfa

$\hspace{0.01em}$


\label{Fact Factor-2.1}

If $ \xda,$ $ \xdb $ are two partitions of $U,$ $ \xda $ a factorization
of $ \xdx,$
and $ \xda $ a refinement of $ \xdb,$ then $ \xdb $ is also a
factorization of $ \xdx.$

\efa

\subparagraph{
Proof
}

$\hspace{0.01em}$


Trivial by definition. $ \xcz $
\\[3ex]

\bfa

$\hspace{0.01em}$


\label{Fact Factor-2.2}

Let $ \xda $ be a factorization of $ \xdx $ over $U,$ $U' \xcc U.$ Then $
\xda \xex U' $ is a factorization
of $ \xdx \xex U' $ over $U'.$

\efa

\subparagraph{
Proof
}

$\hspace{0.01em}$


If $A_{i} \xcs U' \xEd \xCQ,$ let $ \xbs'_{i} \xbe \xdx \xex (A_{i} \xcs
U').$ Let then $ \xbs_{i}:= \xbs'^{+}_{i} \xex A_{i}.$
If $A_{i} \xcs U' = \xCQ,$ let $ \xbs_{i}:= \xbt \xex A_{i}$ for any $
\xbt \xbe \xdx.$
Then $ \xbs:= \xcV \{ \xbs_{i}:i \xbe I\} \xbe \xdx $ by hypothesis, so $
\xbs \xex U' \xbe \xdx \xex U',$ and $ \xbs \xex (A_{i} \xcs U' \}= \xbs
'_{i}.$
$ \xcz $
\\[3ex]

\bfa

$\hspace{0.01em}$


\label{Fact Factor-2.3}

If $\{A,A' \}$ is a factorization of $ \xdx $ over $U,$ $ \xda $ a
factorization of $ \xdx \xex A$ over $ \xCB,$
$ \xda' $ a factorization of $ \xdx \xex A' $ over $A',$ then $ \xda
\xcv \xda' $ is a factorization of $ \xdx $
over $U.$

\efa

\subparagraph{
Proof
}

$\hspace{0.01em}$


Trivial $ \xcz.$
\\[3ex]

\bfa

$\hspace{0.01em}$


\label{Fact Factor-2.4}

If $ \xda,$ $ \xdb $ are two factorizations of $ \xdx,$ then there is a
common
refining factorization.

\efa

\subparagraph{
Proof
}

$\hspace{0.01em}$


Let $ \xbs $ s.t. $ \xcA i \xbe I \xcA j \xbe J(\xbs \xex (A_{i} \xcs
B_{j}) \xbe \xdx \xex (A_{i} \xcs B_{j})),$ show $ \xbs \xbe \xdx.$
Fix $i \xbe I.$ By Fact \ref{Fact Factor-2.2} (page \pageref{Fact Factor-2.2}),
$ \xdb \xex A_{i}$ is a factorization of $ \xdx \xex A_{i},$ so
$ \xcv \{ \xbs \xex (A_{i} \xcs B_{j}):j \xbe J,$ $A_{i} \xcs B_{j} \xEd
\xCQ \}$ $=$ $ \xbs \xex A_{i} \xbe \xdx \xex A_{i}.$ As $ \xda $ is a
factorization
of $ \xdx,$ $ \xbs \xbe \xdx.$ $ \xcz $
\\[3ex]

This does not generalize to infinitely many factorizations:

\be

$\hspace{0.01em}$


\label{Example Factor-2.1}

Take as index set $ \xbo +1,$ all $Y_{k}:=\{0,1\}.$
Take $ \xdx:=\{ \xbs:$ $ \xbs \xex \xbo $ arbitrary, and $ \xbs (\xbo
):=0$ iff $ \xbs \xex \xbo $ is finally constant $\}.$
Consider the partitions $ \xda_{n}:=\{n,(\xbo +1)-n\},$ they are all
fatorizations of $ \xdx,$ as
it suffices to know the sequence from $n+1$ on to know its value on $ \xbo
.$
A common refinement $ \xda $ will have some $A \xbe \xda $ s.t. $ \xbo
\xbe A.$ Suppose there is some
$n \xbe \xbo \xcs A,$ then $A \xcC n+1,$ $A \xcC (\xbo +1)-(n+1),$ this
is impossible, so $A=\{ \xbo \}.$ If $ \xda $
were a factorization of $ \xdx,$ so would be $\{ \xbo,\{ \xbo \}\}$
by Fact \ref{Fact Factor-2.1} (page \pageref{Fact Factor-2.1}), but $ \xdx $
does not factor into $ \xdx \xex \xbo $ and $ \xdx \xex \{ \xbo \}.$

\ee

\bcom

$\hspace{0.01em}$


\label{Comment Factor-2.1}

Above set $ \xdx $ is not definable as a model set of a corresponding
language $ \xdl:$
If $ \xbf $ is not a tautology, there is a model $m$
s.t. $m \xcm \xCN \xbf.$ $ \xbf $ is finite, let its variables be among
$p_{1}, \Xl,p_{n}$ and perhaps
$p_{ \xbo.}$ If $p_{ \xbo }$ is not among its variables, it is trivially
also false in some $m' $
in $ \xdx.$ If it is, then modify $m$ accordingly beyond $n.$ Thus,
exactly all
tautologies are true in $ \xdx,$ but $ \xdx \xEd \xdy =$ the set of all $
\xdl -$models.

\ecom

We have, however:

\bfa

$\hspace{0.01em}$


\label{Fact Factor-2.5}

Let $ \xdx = \xcS \{ \xdx_{m}:m \xbe M\}$ and $ \xdx, \xdx_{m} \xcc \xdy
$ for all $m \xbe M.$

Let $ \xda $ be a partition of $U,$ and a factorization of all $
\xdx_{m}.$

Then $ \xda $ is also a factorization of $ \xdx.$

\efa

\subparagraph{
Proof
}

$\hspace{0.01em}$


Let $ \xbs $ s.t. $ \xcA i \xbe I$ $ \xbs \xex A_{i} \xbe \xdx \xex
A_{i}.$

But $ \xdx \xex A_{i}$ $=$ $(\xcS \{ \xdx_{m}:m \xbe M\}) \xex A_{i}$ $
\xcc $ $ \xcS \{ \xdx_{m} \xex A_{i}:m \xbe M\}:$
Let $ \xbt \xbe \xdx \xex A_{i},$ so by $ \xdx = \xcS \{ \xdx_{m}:m \xbe
M\}$
$ \xbt^{+} \xbe \xdx_{m}$ for all $m \xbe M,$ so $ \xbt \xbe \xdx_{m} \xex
A_{i}$ for all $m \xbe M.$

Thus, $ \xcA i \xbe I, \xcA m \xbe M:$ $ \xbs \xex A_{i} \xbe \xdx_{m}
\xex A_{i},$ so $ \xcA m \xbe M. \xbs \xbe \xdx_{m}$ by prerequisite, so
$ \xbs \xbe \xdx.$ $ \xcz $
\\[3ex]

\bfa

$\hspace{0.01em}$


\label{Fact Factor-2.6}

Let $A \xcv A' $ be a partition of $U,$ and for all $ \xbs \xbe \xdx \xex
A$ and all
$ \xbt:A' \xcp \xcV \{X_{k}:k \xbe A' \}$ with $ \xbt (k) \xbe X_{k}$ $
\xbs \xcv \xbt \xbe \xdx.$ Then

(1) $A \xcv A' $ is a factorization of $ \xdx $ over $U.$

(2) Any partition $ \xda' =\{A'_{k}:k \xbe I' \}$ of $A' $ is a
factorization of $ \xdx \xex A' $ over $A'.$

(3) If $ \xda $ is a factorization of $ \xdx \xex A$ over $ \xCB,$ and $
\xda' $ a partition of $A',$
then $ \xda \xcv \xda' $ is a factorization of $ \xdx.$

\efa

\subparagraph{
Proof
}

$\hspace{0.01em}$


(1) and (2) are trivial, (3) follows from (1), (2), and
Fact \ref{Fact Factor-2.3} (page \pageref{Fact Factor-2.3}). $ \xcz $
\\[3ex]

\bco

$\hspace{0.01em}$


\label{Corollary Factor-2.7}

Let $U=v(\xdl)$ for some language $ \xdl.$ Let $ \xdx $ be definable,
and $\{ \xda_{m}:m \xbe M\}$ be a
set of factorizations of $ \xdx $ over $U.$ Then $ \xda:= \xcv \{
\xda_{m}:m \xbe M\}$ is also a
factorization of $ \xdx.$

\eco

\subparagraph{
Proof
}

$\hspace{0.01em}$


Let $ \xdx =M(T).$ Consider $ \xbf \xbe T.$ $v(\xbf)$ is finite,
consider $ \xdx \xex v(\xbf).$
There are only finitely many different ways $v(\xbf)$ is partitioned by
the
$ \xda_{m},$ let them all be among $ \xda_{m_{0}}, \Xl, \xda_{m_{p}}.$
$M(\xbf) \xex v(\xbf)$ might not be factorized
by all $ \xda_{m_{0}} \xex v(\xbf), \Xl, \xda_{m_{p}} \xex v(\xbf),$
but $M(T) \xex v(\xbf)$ is by
Fact \ref{Fact Factor-2.2} (page \pageref{Fact Factor-2.2}).
By Fact \ref{Fact Factor-2.4} (page \pageref{Fact Factor-2.4}),
$ \xda \xex v(\xbf)$ is a factorization of $M(T) \xex v(\xbf).$

Consider now $ \xdx_{ \xbf }:=(M(T) \xex v(\xbf)) \xCK \xbP \{(0,1):k
\xbe v(\xdl)-v(\xbf)\}.$

By Fact \ref{Fact Factor-2.6} (page \pageref{Fact Factor-2.6}), (1)
$\{v(\xbf),v(\xdl)-v(\xbf)\}$ is a factorization of $ \xdx_{ \xbf }$
over $v(\xdl).$

By Fact \ref{Fact Factor-2.6} (page \pageref{Fact Factor-2.6}), (2)
$ \xda \xex (v(\xdl)-v(\xbf))$ is a factorization of $ \xdx_{ \xbf }
\xex (v(\xdl)-v(\xbf))$
over $v(\xdl)-v(\xbf).$

By Fact \ref{Fact Factor-2.6} (page \pageref{Fact Factor-2.6}), (3)
$ \xda $ is a factorization of $ \xdx_{ \xbf }$ over $v(\xdl).$

$M(T)= \xcS \{(M(T) \xex v(\xbf)) \xCK \xbP \{(0,1):k \xbe v(\xdl)-v(
\xbf)\}$: $ \xbf \xbe T\},$ so
by Fact \ref{Fact Factor-2.5} (page \pageref{Fact Factor-2.5}),
$ \xda $ is a factorization of $M(T).$

$ \xcz $
\\[3ex]

\bcom

$\hspace{0.01em}$


\label{Comment Factor-2.2}

Obviously, it is unimportant here that we have only 2 truth values, the
proof
would just as well work with any, even an infinite, number of truth
values.
What we really need is the fact that a formula affects only finitely many
propositional variables, and the rest are free.

\ecom

Unfortunately, the manner of coding can determine if there is a
factorization,
as can be seen by the following example:

\be

$\hspace{0.01em}$


\label{Example Factor-2.2}

 \xEh

 \xDH
p= ``$blue$'', q= ``$round$'', q'= ``$blue$ $iff$ $round$''.

Then

$p \xcu q$ = $blue$ $and$ $round$,  $\neg p \xcu \neg q$ =
$\neg blue$ $and$  $\neg round$

$p \xcu q'$ = $blue$ $and$ $round$,  $\neg p \xcu q'$ =
$ \neg blue$ $and$ $\neg round$

Thus, both code the same (meta-) situation, the first cannot be
factorized,
the second can.

Our example (first presented in  \cite{Sch07a})
is discussed in more detail in
 \cite{Mak09}, see Section 5. there.

 \xDH

More generally, we can code e.g. the non-factorising situation
$\{p \xcu q \xcu r, \xCN p \xcu \xCN q \xcu \xCN r\}$ also using $q' =p
\xcr q,$ $r' =p \xcr r,$ and have then the
factorising situation $\{p \xcu q \xcu r, \xCN p \xcu q' \xcu r' \}.$

 \xDH

The following situation cannot be made factorising:
$\{p \xcu q,$ $p \xcu \xCN q,$ $ \xCN p \xcu \xCN q\}.$ Suppose there were
some such solution.
Then we need some $p' $ and $q',$ and all 4 possibilities
$\{p' \xcu q',$ $p' \xcu \xCN q',$ $ \xCN p' \xcu q',$ $ \xCN p' \xcu
\xCN q' \}.$ If we do not admit impossible
situations (i.e. one of the 4 possibilities is a contradictory coding),
then
2 possibilities have to contain the same situation, e.g. $p \xcu q.$ But
they
are mutually exclusive (as they are negations), so this is impossible.

 \xEj

\ee

$ \xcz $
\\[3ex]

\br

$\hspace{0.01em}$


\label{Remark Generalisations}

As we worked with abstract sequences, which need not be models,
we can apply our results and the ideas behind them (essentially
due to Parikh/Rodrigues) e.g., to

 \xEI

 \xDH

Update: if a set of sequences (where the points are now models, and not
true/false) factorizes, then we can update the components (i.e. look for
the locally ``best'' subsequences), and then compose them to the
globally ``best'' sequences.

 \xDH

Utility streams: if a set of utility streams factorizes, we can do the
same, commutativity and associativity of addition will guarantee the
desired result.

 \xDH

Preferential reasoning: Again, we factorize, and choose the locally
best which we compose to the globally best.

 \xEJ

The idea is always the same: If a set factorizes, choose locally, and
compose to the global choice - provided this is the desired result!
\section{
Factorisation and Hamming distance
}

\er

$ \xCO $
\index{Hamming-Definition}
\label{Hamming-Definition}

\bd

$\hspace{0.01em}$


\label{Definition Hamming-Distance}

Given $x,y \xbe \xbS,$ a set of sequences over an index set $I,$ the
Hamming distance
comes in two flavours:

$d_{s}(x,y):=\{i \xbe I:x(i) \xEd y(i)\},$ the set variant,

$d_{c}(x,y):=card(d_{s}(x,y)),$ the counting variant.

We define $d_{s}(x,y) \xck d_{s}(x',y')$ iff $d_{s}(x,y) \xcc d_{s}(x'
,y'),$

thus, $s-$distances are not always comparabel.

We can also give different importance to different $i$ in the counting
variant, so e.g.,
$d_{c}(\xBc x,x'  \xBe, \xBc y,y'  \xBe)$ might be 1 if $x \xEd y$ and $x'
=y',$ but 2 if
$x=y$ and $x' \xEd y'.$

\ed

$ \xCO $

\bfa

$\hspace{0.01em}$


\label{Fact Hamming-Distance}

$d_{c}$ has the normal addition, set union takes the role of addition for
$d_{s},$
$ \xCQ $ takes the role of 0 for $d_{s},$
both are distances in the following sense:

(1) $d(x,y)=0$ iff $x=y,$

(2) $d(x,y)=d(y,x),$

(3) the triangle inequality holds, for the set variant in the form
$d_{s}(x,z) \xcc d_{s}(x,y) \xcv d_{s}(y,z).$

\efa

\subparagraph{
Proof
}

$\hspace{0.01em}$


(3) If $i \xce d_{s}(x,y) \xcv d_{s}(y,z),$ then $x(i)=y(i)=z(i),$ so
$x(i)=z(i)$ and
$i \xce d_{s}(x,z).$

The others are trivial.

$ \xcz $
\\[3ex]

Both Hamming distances cooperate well with factorization, as we will see
now. This is not surprising, as Hamming distances work componentwise.

\bd

$\hspace{0.01em}$


\label{Definition Parikh-General}

We say that a revision function $*$ factorizes iff for all $T$ and $ \xbf
$ and
joint factorisations, which we write for simplicity (and immediately for
models)
$M(T)=M(T) \xex \xdl_{1} \xCK, \Xl, \xCK M(T) \xex \xdl_{n},$ $M(\xbf
)=M(\xbf) \xex \xdl_{1} \xCK, \Xl, \xCK M(\xbf) \xex \xdl_{n},$
$M(T* \xbf)$ $=$ $M(T) \xfA M(\xbf)$ $=$
$((M(T) \xex \xdl_{1}) \xfA (M(\xbf) \xex \xdl_{1})) \xCK, \Xl, \xCK
((M(T) \xex \xdl_{n}) \xfA (M(\xbf) \xex \xdl_{n})).$

To simplify notation, we will speak about $ \xbS \xfA T,$ $ \xbS_{i} \xfA
T_{i},$
$ \xbs_{1} \xCK  \Xl  \xbs_{n},$ etc.

\ed

The advantage is that, when factorisation is possible, we can work with
smaller theories, formulas, and languages, and then do a trivial
composition operation by considering the product.

\bfa

$\hspace{0.01em}$


\label{Fact Parikh-Hamming}

If $*$ is defined by the counting or the set variant of the Hamming
distance,
then $*$ factorises.

\efa

\subparagraph{
Proof
}

$\hspace{0.01em}$


We do the proof for the set variant, the counting variant
proof is similar.

Let a factorisation as in the definition be given, and suppose
$ \xbt \xbe T$ has minimal distance from $ \xbS,$ i.e. $ \xbt \xbe \xbS
\xfA T.$
We show that each $ \xbt_{i}$ has minimal
distance from $ \xbS_{i}.$ If not, there is $ \xbt'_{i}$ closer to $
\xbS_{i},$ but then
$ \xbt',$ which is like $ \xbt,$ only $ \xbt_{i}$ is replaced by $ \xbt
'_{i}$ is also in $T,$
by factorisation. By definition of the Hamming distance, $ \xbt' $ is
closer
to $ \xbS $ than $ \xbt $ is, $contradiction.$ Thus $ \xbt \xbe (\xbS_{1}
\xfA T_{1}) \xCK  \Xl  \xCK (\xbS_{n} \xfA T_{n}).$
Conversely, let all $ \xbt_{i} \xbe (\xbS_{i} \xfA T_{i}),$ we have to
show that
$ \xbt:= \xbt_{1} \xCK  \Xl  \xCK \xbt_{n} \xbe \xbS \xfA T.$ By
factorisation, $ \xbt \xbe T.$ If there were a
closer $ \xbt' \xbe T,$ then at least one of the components $ \xbt'_{i}$
would be closer
than $ \xbt_{i},$ $contradiction.$

$ \xcz $
\\[3ex]

The authors do not know if all factorising distance defined revisions can
be
defined by one of the above Hamming distances.
\section{
Preferential modelling of defaults
}

Reiter defaults have the advantage to give results also for non-ideal
cases. If, by default, $ \xba $ and $ \xba' $ hold, but $ \xba $ is
inconsistent with the
current situation, then $ \xba' $ will still ``fire''. Preferential
structures
say nothing about non-ideal cases. We construct special preferential
structures
which have the same behaviour as Reiter defaults. In addition, specificity
will be used to solve conflicts.

The idea is simple.

For simplicity, we admit direct contradictions: $ \xbf \xcn \xbq $ and $
\xbf \xcn \xCN \xbq.$ This
is done only to make the representation proof simple. One can do without,
but
pays with more complexity (see below). We also use structures with one
copy
of each model only.

\bd

$\hspace{0.01em}$


\label{Definition Default-Specificity}

(1) We call the default $ \xbf \xcn \xbq $ more specific than the default
$ \xbf' \xcn \xbq' $ iff
$ \xbf \xcl \xbf' $ - $ \xcl $ is classical consequence.

(2) We say that the default $ \xbf \xcn \xbq $ separates $m$ and $m' $ iff
$m,m' \xcm \xbf,$ but only
one of $m,m' $ satisfies $ \xbq.$

\ed

Consider two models, $m,m'.$ Take the most specific defaults which
separate them. For each such default $ \xbf \xcn \xbq,$ if $m \xcm \xbq
,$ $m' \xcM \xbq,$ then set
$m \xeb m'.$ We might introduce cycles of length 2 here.

\br

$\hspace{0.01em}$


\label{Remark Ind-1}

(1) The construction has a flavour of rankedness, as, if possible,
we make each ``good'' element smaller than each
``less good'' element. If e.g. $ \xbf \xcn \xbq $ is the default, $m,m'
,n,n' \xcm \xbf,$ and $m,n \xcm \xbq,$
but $m',n' \xcM \xbq,$ then $m \xeb m',$ $n \xeb m',$ $m \xeb n',$ $n
\xeb n'.$

(2) We may create indirectly loops, as we may have $m \xeb m' \xeb m'',$
and also see
$m'' \xeb m.$

(3) Let $m \xbe X$ be minimal in our construction. Then there is no
default
$ \xbf \xcn \xbq $ s.t. $m \xcm \xbf \xcu \xCN \xbq,$ and there is $m'
\xbe X,$ $m' \xcm \xbf \xcu \xbq.$ Thus, minimal elements
are ``as good as possible'', i.e. there is no better one in $X.$ Of course,
there might be a default $ \xbf \xcn \xbq $ with $m \xcm \xbf \xcu \xCN
\xbq,$ but
``better'' elements are outside $X.$ Thus, minimal elements are an
approximation of the ideal case. We can also consider minimal elements
as a revision of the ideal case by $X,$ in the sense that we cannot get
closer to the ideal within $X.$

\er

Conversely, given $ \xCf any$ preferential structure, we take any two
models $m,m',$
if $m \xeb m' $ (which we see by $m' \xce \xbm (\{m,m' \})),$ we add the
default $Th\{m,m' \} \xcn Th\{m\}.$
By the basic law of 1-copy preferential structures, we create the
structure
again. (Note that $\{m,m' \}$ is the most specific set containing both.)
Thus, our approach cannot result in new structural rules for preferential
structures, like smoothness, rankedness, etc.
\section{
Remarks on independence
}

The idea of independence was realized for defaults by trying to satisfy
them independently, so if one fails, the others still have a chance.
This is a simple idea.

The case of theory revision is more complicated, as we have no predefined
structure. In particular, the starting theory $T$ can be just a
``blob'' which makes independence difficult to realize. Moreover, we
may try to revise just once, so no multiple default satisfaction or so
is needed. Perhaps the Parikh idea, and its refinement through (modified)
Hamming distances is all one can achieve.

It is evident how to treat our form of independent revision with IBRS, it
can just be written as a diagram.

Perhaps the best way to write defaults as a diagram, is the trivial one:
$ \xba \xcn \xbb $ will be written $ \xba \xch \xbb,$ and the treatment
is in the evaluation of
the diagram - as outlined above.
\subsection{
Epistemic states and independence
}

It is probably adequate to say that  \cite{AGM85} consider the
revision
function $*$ an epistemic state (depending on $K),$ as revealed by the
notion of epistemic entrenchment, and its equivalence to a revision
function.
In  \cite{LMS01}, the global distance can probably be seen as
(fixed, global)
epistemic state.
Essentially the critique of such too rigid, fixed, epistemic states
resulted
in dynamic states of  \cite{DP94} and  \cite{DP97}.
The approach in  \cite{Spo88} incorporated already a dynamic
approach.
An excellent short overview of such dynamic revision approaches can be
found
in  \cite{Ker99}.

In preferential structures, we might see the relation choosing the normal
situations as (again fixed) epistemic state. In counterfactual
conditionals,
the distance can again be seen as the underlying epistemic state.

The present authors see the higher order arrows of reactive structures
(see e.g.,  \cite{GS08b}) as expressing epistemic states or changes
of
epistemic states. This will be explored in future research by the
present authors.

But, we can also see the approaches discussed in this article as
an epistemic state, which can perhaps be resumed as:
``divide and conquer''.

\vspace{3mm}



\vspace{3mm}




\begin{thebibliography}{xxxxxx}

\addcontentsline{toc}{section}{References}


\bibitem[AGM85]{AGM85}
C.Alchourron, P.Gardenfors, D.Makinson,
``On the Logic of Theory Change: partial meet contraction and
revision functions'', Journal of Symbolic Logic, Vol. 50,
pp. 510-530, 1985

\bibitem[CP00]{CP00}
S.Chopra, R.Parikh, ``Relevance sensitive belief structures'',
Annals of Mathematics and Artificial Intelligence,
vol. 28, No. 1-4, pp. 259-285, 2000

\bibitem[DP94]{DP94}
A.Darwiche, J.Pearl, ``On the Logic of Iterated Belief Revision'',
in: ``Proceedings of the fifth Conference on Theoretical Aspects of
Reasoning about Knowledge'', R.Fagin ed., pp. 5-23, Morgan Kaufman,
Pacific Grove, CA, 1994

\bibitem[DP97]{DP97}
A.Darwiche, J.Pearl, ``On the Logic of Iterated Belief Revision'',
Journal of Artificial Intelligence, Vol. 89, No. 1-2, pp. 1-29,
1997

\bibitem[GS08b]{GS08b}
D.Gabbay, K.Schlechta, ``Reactive preferential structures and
nonmonotonic consequence'',
to appear in Review of Symbolic Logic
hal-00311940, arXiv 0808.3075

\bibitem[Ker99]{Ker99}
G.Kern-Isberner, ``Postulates for conditional belief revision'',
Proceedings IJCAI 99, T.Dean ed., Morgan Kaufmann, pp.186-191, 1999

\bibitem[LMS01]{LMS01}
D.Lehmann, M.Magidor, K.Schlechta: ``Distance Semantics for Belief
Revision'', Journal of Symbolic Logic, Vol.66, No. 1, March 2001, $p.$
295-317

\bibitem[Mak09]{Mak09}
D.Makinson: ``Propositional relevance through letter-sharing''
to appear in Journal of Applied Logic, special issue, J.Delgrande ed.

\bibitem[Rod97]{Rod97}
O.T.Rodrigues: ``A methodology for iterated information change'',
PhD thesis, Imperial College, London, 1997

\bibitem[Sch07a]{Sch07a}
K.Schlechta: ``Factorization'',
HAL, arXiv.org 0712.4360v1, 2007

\bibitem[Spo88]{Spo88}
W.Spohn, ``Ordinal conditional functions: A dynamic theory of
epistemic states''. In: W.L.Harper and B.Skyrms, (eds.), ``Causation in
Decision, Belief Change, and Statistics'', vol. 2, p.105-134, Reidel,
Dordrecht 1988,

\end{thebibliography}
\end{document}